\documentclass[12pt]{article}%
\usepackage{amsfonts}
\usepackage{amssymb}
\usepackage{amsmath}
\usepackage{amscd}
\usepackage{latexsym}
\usepackage[onehalfspacing]{setspace}
\usepackage{theorem}
\usepackage{scalefnt}
\usepackage{ifthen}
\usepackage{natbib}
\usepackage[colorlinks=true,urlcolor=black,linkcolor=black,pdfpagemode="None",pdfstartview=FitH]%
{hyperref}
\usepackage{graphicx}
\usepackage[font={footnotesize},textfont={footnotesize}]{caption}
\usepackage{subfig}%
\providecommand{\U}[1]{\protect\rule{.1in}{.1in}}
\newtheorem{theorem}{Theorem}

\newtheorem{lemma}{Lemma}

\newcommand{\I}{{\rm 1\hspace*{-0.4ex}\rule{0.1ex}{1.52ex}\hspace*{0.2ex}}}
\hypersetup{pdfborder = {0 0 0},colorlinks=true,linkcolor=blue,citecolor=blue}

\renewcommand{\cite}{\citet*}
\begin{document}

\title{\vspace{-0.5in}Alternative Asymptotics and the Partially Linear Model with
Many Regressors\thanks{The authors thank comments from Alfonso Flores-Lagunes,
Lutz Kilian, seminar participants at Bristol, Brown, Cambridge, Exeter,
Indiana, LSE, Michigan, MSU, NYU, Princeton, Rutgers, Stanford, UCL, UCLA,
UCSD, UC-Irvine, USC, Warwick and Yale, and conference participants at the
2010 Joint Statistical Meetings and the 2010 LACEA Impact Evaluation Network
Conference. The first author gratefully acknowledges financial support from
the National Science Foundation (SES 1122994). The second author gratefully
acknowledges financial support from the National Science Foundation (SES
1124174) and the research support of CREATES (funded by the Danish National
Research Foundation). The third author gratefully acknowledges financial
support from the National Science Foundation (SES 1132399).}\bigskip}
\author{Matias D. Cattaneo\thanks{Department of Economics, University of Michigan.}
\and Michael Jansson\thanks{Department of Economics, UC Berkeley and \emph{CREATES}%
.}
\and Whitney K. Newey\thanks{Department of Economics, MIT.}}
\maketitle

\begin{abstract}
Non-standard distributional approximations have received considerable
attention in recent years. They often provide more accurate approximations in
small samples, and theoretical improvements in some cases. This paper shows
that the seemingly unrelated \textquotedblleft many
instruments\ asymptotics\textquotedblright\ and \textquotedblleft small
bandwidth\ asymptotics\textquotedblright\ share a common structure, where the
object determining the limiting distribution is a V-statistic with a remainder
that is an asymptotically normal degenerate U-statistic. We illustrate how
this general structure can be used to derive new results by obtaining a new
asymptotic distribution of a series estimator of the partially linear model
when the number of terms in the series approximation possibly grows as fast as
the sample size, which we call \textquotedblleft many
terms\ asymptotics\textquotedblright.

\end{abstract}

\textbf{JEL classification}\textit{:} C13, C31.\medskip

\textbf{Keywords}: non-standard asymptotics, partially linear model, many
terms, adjusted variance.\bigskip%

\thispagestyle{empty}
\setcounter{page}{0}
\newpage

\section{Introduction\label{section:intro}}

Many instrument asymptotics, where the number of instruments grows as fast as
the sample size, has proven useful for instrumental variables
(IV)\ estimators. \cite{Kunitomo_1980_JASA} and \cite{Morimune_1983_ECMA}
derived asymptotic variances that are larger than the usual formulae when the
number of instruments and sample size grow at the same rate, and
\cite{Bekker_1994_ECMA} and others provided consistent estimators of these
larger variances. \cite{Hansen-Hausman-Newey_2008_JBES} showed that using many
instrument standard errors provides a theoretical improvement for a range of
number of instruments and a practical improvement for estimating the returns
to schooling. Thus, many instrument asymptotics and the associated standard
errors have been demonstrated to be a useful alternative to the usual
asymptotics for instrumental variables.

Instrumental variable estimators implicitly depend on a nonparametric series
estimator. Many instrument asymptotics has the number of series terms growing
so fast that the series estimator is not consistent. Analogous asymptotics for
kernel-based density-weighted average derivative estimators has been
considered by
\citet{Cattaneo-Crump-Jansson_2010_JASA,Cattaneo-Crump-Jansson_2014a_ET}%
. They show that when the bandwidth shrinks faster than needed for consistency
of the kernel estimator, the variance of the estimator is larger than the
usual formula. They also find that correcting the variance provides an
improvement over standard asymptotics for a range of bandwidths.

The purpose of this paper is to show that these results share a common
structure, and to illustrate how this structure can be used to derive new
results. The common structure is that the object determining the limiting
distribution is a V-statistic, which can be decomposed into a bias term, a
sample average, and a \textquotedblleft remainder\textquotedblright\ that is
an asymptotically normal degenerate U-statistic. Asymptotic normality of the
remainder distinguishes this setting from other ones involving V-statistics.
Here the asymptotically normal remainder comes from the number of series terms
going to infinity or bandwidth shrinking to zero, while the behavior of a
degenerate U-statistic tends to be more complicated in other settings. When
the number of terms grows as fast as the sample size, or the bandwidth shrinks
to zero at an appropriate rate, the remainder has the same magnitude as the
leading term, resulting in an asymptotic variance larger than just the
variance of the leading term. The many instrument and small bandwidth results
share this structure. In keeping with this common structure, we will
henceforth refer to such results under the general heading of
\textquotedblleft alternative asymptotics\textquotedblright.

The alternative asymptotics that we discuss in this paper applies to
statistics that take a specific V-statistic representation, or may be
approximated by it sufficiently accurately, and therefore it does not apply
broadly to all possible semiparametric settings. Nonetheless, as we illustrate
below, this structure arises naturally in several interesting problems in
Economics and Statistics. In particular, we show formally that applying this
common structure to a series estimator of the partially linear model leads to
new results. These results allow the number of terms in the series
approximation to grow as fast as the sample size. The asymptotic distribution
of the estimator is derived and it is shown to have a larger asymptotic
variance than the usual formula, which is in fact a natural and generic
consequence of the specific structure that we highlight in this paper. We also
find that under homoskedasticity, the classical degrees of freedom adjusted
homoskedastic standard error estimator from linear models is consistent even
when the number of terms is \textquotedblleft large\textquotedblright%
\ relative to the sample size. This result offers a large sample, distribution
free justification for the degrees of freedom correction when many series
terms are employed. Constructing automatic consistent standard error estimator
under (conditional) heteroskedasticity of unknown form in this setting turns
out to be quite challenging. In \cite{Cattaneo-Jansson-Newey_2015_HCSEManyCov}%
, we present a detailed discussion of heteroskedasticity-robust standard
errors for general linear models with increasing dimension, which covers the
partially linear model with many terms studied herein as a special case.

The rest of the paper is organized as follows. Section
\ref{section:CommonStructure} describes the common structure of many
instrument and small bandwidth asymptotics, and also shows how the structure
leads to new results for the partially linear model. Section
\ref{section:ManyReg} formalizes the new distributional approximation for the
partially linear model. Section \ref{section:simuls} reports results from a
small simulation study aimed to illustrate our results in small samples.
Section \ref{section:conclusion} concludes. The appendix collects the proofs
of our results.

\section{A Common Structure\label{section:CommonStructure}}

To describe the common structure of many instrument and small bandwidth
asymptotics, let $W_{1},\ldots,W_{n}$ denote independent random vectors. We
consider an estimator $\hat{\beta}$ of a generic parameter of interest
$\beta_{0}\in\mathbb{R}^{d}$ satisfying%
\begin{equation}
\sqrt{n}(\hat{\beta}-\beta_{0})=\hat{\Gamma}_{n}^{-1}S_{n},\text{\qquad}%
S_{n}=\sum_{1\leq i,j\leq n}u_{ij}^{n}(W_{i},W_{j}),\label{eq:vstat}%
\end{equation}
where $u_{ij}^{n}(\cdot)$ is a function that can depend on $i$, $j$, and $n$.
We allow $u$ to depend on $n$ to account for number of terms or bandwidths
that change with the sample size. Also, we allow $u$ to vary with $i$ and $j$
to account for dependence on variables that are being conditioned on in the
asymptotics, and so treated as nonrandom.

We assume throughout this section that there exists a sequence of non-random
matrices $\Gamma_{n}$ satisfying $\Gamma_{n}^{-1}\hat{\Gamma}_{n}%
\rightarrow_{p}I_{d}$ for $I_{d}$ the $d\times d$ identity matrix, and hence
we focus on the V-statistic $S_{n}$. (All limits are taken as $n\rightarrow
\infty$ unless explicitly stated otherwise.)\ This V-statistic has a well
known (Hoeffding-type) decomposition that we describe here because it is an
essential feature of the common structure. For notational implicitly we will
drop the $W_{i}$ and $W_{j}$ arguments and set $u_{ij}^{n}=u_{ij}^{n}%
(W_{i},W_{j})$ and $\tilde{u}_{ij}^{n}=u_{ij}^{n}+u_{ji}^{n}-\mathbb{E}%
[u_{ij}^{n}+u_{ji}^{n}]$.

Letting $\Vert\cdot\Vert$ denote the Euclidean norm, and if $\mathbb{E}[\Vert
u_{ij}^{n}\Vert]<\infty$\ for all $i,j,n$, then%
\begin{equation}
S_{n}=B_{n}+\Psi_{n}+U_{n},\label{eq:decomp}%
\end{equation}
where%
\[
B_{n}=\mathbb{E}[S_{n}],\text{\qquad}\Psi_{n}=\sum_{1\leq i\leq n}\psi_{i}%
^{n}(W_{i}),\text{\qquad}U_{n}=\sum_{2\leq i\leq n}D_{i}^{n}(W_{i},...,W_{1}),
\]%
\[
\psi_{i}^{n}(W_{i})=u_{ii}^{n}-\mathbb{E}[u_{ii}^{n}]+\sum_{1\leq j\leq
n,j\neq i}\mathbb{E}[\tilde{u}_{ij}^{n}|W_{i}],
\]%
\[
D_{i}^{n}(W_{i},...,W_{1})=\sum_{1\leq j\leq n,j<i}(\tilde{u}_{ij}%
^{n}-\mathbb{E}[\tilde{u}_{ij}^{n}|W_{i}]-\mathbb{E}[\tilde{u}_{ij}^{n}%
|W_{j}]),
\]
so that $\mathbb{E}[\psi_{i}^{n}(W_{i})]=0$, $\mathbb{E}[D_{i}^{n}%
(W_{i},...,W_{1})|W_{i-1},...,W_{1}]=0$, and $\mathbb{E}[\Psi_{n}U_{n}]=0$.
This decomposition of a V-statistic is well known (e.g., \cite[Chapter
11]{vanderVaart_1998_Book}), and shows that $S_{n}$ can be decomposed into a
sum $\Psi_{n}$ of independent terms, a U-statistic remainder $U_{n}$ that is a
martingale difference sum and uncorrelated with $\Psi_{n}$, and a pure bias
term $B_{n}$.\footnote{In time series contexts, the exact decomposition is
less useful, but approximations thereof with properties similar to those we
discuss herein can be developed. For an example and related references see
\cite{Atchade-Cattaneo_2014_SPA}.} The decomposition is important in many of
the proofs of asymptotic normality of semiparametric estimators, including
\cite{Powell-Stock-Stoker_1989_ECMA}, with the limiting distribution being
determined by $\Psi_{n}$, and $U_{n}$ being treated as a \textquotedblleft
remainder\textquotedblright\ that is of smaller order under a particular
restriction on the tuning parameter sequence (e.g., when the bandwidth shrinks
slowly enough).

An interesting feature of the decomposition (\ref{eq:decomp}) in
semiparametric settings is that $U_{n}$ is asymptotically normal at some rate
when the number of series terms grow or the bandwidth shrinks to zero. To be
specific, under regularity conditions and appropriate tuning parameter
sequences that we make precise below, it turns out that%
\[
\left[
\begin{array}
[c]{c}%
\mathbb{V}[\Psi_{n}]^{-1/2}\Psi_{n}\\
\mathbb{V}[U_{n}]^{-1/2}U_{n}%
\end{array}
\right]  \rightarrow_{d}\mathcal{N}(0,I_{2d}).
\]
In other settings, where the underlying kernel of the U-statistic does not
vary with the sample size, the asymptotic behavior of $U_{n}$ is usually more
complicated: because it is a degenerate U-statistic, it would converge to a
weighted sum of independent chi-square random variables (e.g., \cite[Chapter
12]{vanderVaart_1998_Book}). However, in semiparametric-type settings as those
considered here, the kernel of the underlying U-statistic forming $U_{n}$
changes with the sample size and hence, under particular tuning parameter
configurations, the individual contributions $D_{i}^{n}(W_{i},...,W_{1})$ to
$U_{n}$ can be made small enough to satisfy a Lindeberg-Feller condition and
thus obtain a Gaussian limiting distribution (usually employing the martingale
property of $U_{n}$). For an interesting discussion of this phenomenon, see
\cite{deJong_1987_PTRF}. The asymptotic normality property of $U_{n}$ has been
shown for certain classes of both series and kernel based estimators, as
further explained below.

Alternative asymptotics occurs when the number of series terms grows or the
bandwidth shrinks fast enough so that $\mathbb{V}[\Psi_{n}]$ and
$\mathbb{V}[U_{n}]$ have the same magnitude in the limit. Because of
uncorrelatedness of $\Psi_{n}$ and $U_{n}$, the asymptotic variance will be
larger than the usual formula which is $\lim_{n\rightarrow\infty}%
\mathbb{V}[\Psi_{n}]$ (assuming the limit exists). As a consequence,
consistent variance estimation under alternative asymptotics requires
accounting for the contribution of $U_{n}$ to the (asymptotic)\ sampling
variability of the statistic. Accounting for the presence of $U_{n}$ should
also yield improvements when numbers of series terms and bandwidths do not
satisfy the knife-edge conditions of alternative asymptotics, since $U_{n}$ is
part of the semiparametric statistic. For instance, if the number of series
terms grows just slightly slower than the sample size then accounting for the
presence of $U_{n}$ should still give a better large sample approximation.
\cite{Hansen-Hausman-Newey_2008_JBES} show such an improvement for many
instrument asymptotics. It would be good to consider such improved
approximations more generally, though it is beyond the scope of this paper to
do so.

Distribution theory under alternative asymptotics may be seen as a
generalization of the conventional large sample distributional approximation
approach in the sense that under conventional sequences of tuning parameters
the asymptotic variances emerging from both approaches coincide. But, the
alternative asymptotic approximation also allows for other tuning parameter
sequences and, in this case, the limiting asymptotic variance is seen to be
larger than usual. Thus, in general, there is no reason to expect that the
usual standard error formulas derived under conventional asymptotics will
remain valid more generically. From this perspective, alternative asymptotics
are useful to provide theoretical justification for new standard error
formulas that are consistent under more general sequences of tuning
parameters, that is, under both conventional and alternative asymptotics. We
refer to the latter standard error formulas as being more robust than the
usual standard error formulas available in the literature. For instance, using
these ideas, the need for new, more robust standard errors formulas was made
before for many instrument asymptotics in IV\ models
(\cite{Hansen-Hausman-Newey_2008_JBES}) and small bandwidth asymptotics in
kernel-based semiparametrics (\cite{Cattaneo-Crump-Jansson_2014a_ET}).

To illustrate these ideas, we show next that both many instrument asymptotics
and small bandwidth asymptotics have the structure described above, and we
also employ this approach to derive new results in the case of a series
estimator of the partially linear model, which we refer to as
\textquotedblleft many terms asymptotics\textquotedblright.

\subsection*{Example 1: \textquotedblleft Many Instrument
Asymptotics\textquotedblright}

The first example is concerned with the case of many instrument asymptotics.
For simplicity we focus on the JIVE2 estimator of
\cite{Angrist-Imbens-Krueger_1999_JAE}, but the idea applies to other IV
estimators such as the limited information maximum likelihood estimator. See
\cite{Chao-Swanson-Hausman-Newey-Woutersen_2012_ET} for more details,
including regularity conditions under which the following discussion can be
made rigorous.

Let $(y_{i},x_{i}^{\prime},z_{i}^{\prime})^{\prime}$, $i=1,\ldots,n$, be a
random sample generated by the model%
\begin{equation}
y_{i}=x_{i}^{\prime}\beta_{0}+\varepsilon_{i},\text{\hspace{0.2in}}%
\mathbb{E}[\varepsilon_{i}|z_{i}]=0,\label{eq:IVmodel}%
\end{equation}
where $y_{i}$ is a scalar dependent variable, $x_{i}\in\mathbb{R}^{d}$ is a
vector of endogenous variables, $\varepsilon_{i}$ is a disturbance, and
$z_{i}\in\mathbb{R}^{K}$ is a vector of instrumental variables.

To describe the JIVE2 estimator of $\beta_{0}$ in (\ref{eq:IVmodel}), let
$Q_{ij}$ denote the $(i,j)$-th element of $Q=Z(Z^{\prime}Z)^{-1}Z^{\prime}$,
where $Z=[z_{1},\cdots,z_{n}]^{\prime}$. After centering and scaling, the
JIVE2 estimator $\hat{\beta}$ satisfies%
\[
\sqrt{n}(\hat{\beta}-\beta_{0})=(\frac{1}{n}\sum_{1\leq i,j\leq n,j\neq
i}Q_{ij}x_{i}x_{j}^{\prime})^{-1}(\frac{1}{\sqrt{n}}\sum_{1\leq i,j\leq
n,j\neq i}Q_{ij}x_{i}\varepsilon_{j}).
\]
Conditional on $Z,$ $\hat{\beta}$ has the structure in (\ref{eq:vstat}) with
$W_{i}=(x_{i}^{\prime},\varepsilon_{i})^{\prime}$ and%
\[
\hat{\Gamma}_{n}=\frac{1}{n}\sum_{1\leq i,j\leq n,j\neq i}Q_{ij}x_{i}%
x_{j}^{\prime},\text{\qquad}u_{ij}^{n}(W_{i},W_{j})=%
\I
(i\neq j)Q_{ij}x_{i}\varepsilon_{j}/\sqrt{n},
\]
where $%
\I
(\cdot)$ is the indicator function.

For $i\neq j$, $\mathbb{E}[u_{ij}^{n}(W_{i},W_{j})|Z]=0$ and%
\[
\mathbb{E}[u_{ij}^{n}(W_{i},W_{j})|W_{i},Z]=Q_{ij}x_{i}\mathbb{E}%
[\varepsilon_{j}|Z]=0,\text{\hspace{0.2in}}\mathbb{E}[u_{ji}^{n}(W_{j}%
,W_{i})|W_{i},Z]=Q_{ij}\Upsilon_{j}\varepsilon_{i}/\sqrt{n},
\]
where $\Upsilon_{i}=\mathbb{E}[x_{i}|z_{i}]$ can be interpreted as the reduced
form for observation $i$. As a consequence, (\ref{eq:decomp}) is satisfied
with $B_{n}=0$,%
\[
\psi_{i}^{n}(W_{i})=(\sum_{1\leq j\leq n,j\neq i}Q_{ij}\Upsilon_{j}%
)\varepsilon_{i}=\Upsilon_{i}(1-Q_{ii})\varepsilon_{i}/\sqrt{n}-(\Upsilon
_{i}-\sum_{1\leq j\leq n}Q_{ij}\Upsilon_{j})\varepsilon_{i}/\sqrt{n},
\]%
\[
D_{i}^{n}(W_{i},...,W_{1})=\sum_{1\leq j\leq n,j<i}Q_{ij}\left(
v_{i}\varepsilon_{j}+v_{j}\varepsilon_{i}\right)  /\sqrt{n},\text{\hspace
{0.2in}}v_{i}=x_{i}-\Upsilon_{i}.
\]

Because $\Upsilon_{i}-\sum_{j=1}^{n}Q_{ij}\Upsilon_{j}$ is the $i$-th residual
from regressing the reduced form observations on $Z$, by appropriate
definition of the reduced form this can generally be assumed to vanish as the
sample size grows. In that case,%
\[
\Psi_{n}=\frac{1}{\sqrt{n}}\sum_{1\leq i\leq n}\Upsilon_{i}(1-Q_{ii}%
)\varepsilon_{i}+o_{p}(1).
\]
Furthermore, under standard asymptotics $Q_{ii}$ will go to zero, so the
limiting variance of the leading term in $\Psi_{n}$ corresponds to the usual
asymptotic variance for IV. The degenerate U-statistic term is%
\[
U_{n}=\frac{1}{\sqrt{n}}\sum_{1\leq i,j\leq n,j<i}Q_{ij}\left(  v_{i}%
\varepsilon_{j}+v_{j}\varepsilon_{i}\right)  .
\]
\cite{Chao-Swanson-Hausman-Newey-Woutersen_2012_ET} apply a martingale central
limit theorem to show that this $U_{n}$ will be asymptotically normal when
$K\rightarrow\infty$ and certain regularity conditions hold. The conditions of
the martingale central limit theorem are verified by showing that certain
linear combinations with coefficients depending on the elements of $Q$ go to
zero as $K\rightarrow\infty$. In the proof, this makes individual terms
asymptotically negligible, with a Lindeberg-Feller condition being satisfied.
Alternative asymptotics occurs when $K$ grows as fast as $n$, resulting in
$\mathbb{V}[\Psi_{n}]$ and $\mathbb{V}[U_{n}]$ having the same magnitude in
the limit.

\subsection*{Example 2: \textquotedblleft Small Bandwidth
Asymptotics\textquotedblright}

The second example shows that small bandwidth asymptotics for certain
kernel-based semiparametric estimators also has the structure outlined above.
To keep the exposition simple we focus on an estimator of the integrated
squared density, but the structure of this estimator is shared by the
density-weighted average derivative estimator of
\cite{Powell-Stock-Stoker_1989_ECMA} treated in
\cite{Cattaneo-Crump-Jansson_2014a_ET} and more generally by estimators of
density-weighted averages and ratios thereof (see, e.g., \cite[Section
2]{Newey-Hsieh-Robins_2004_ECMA} and references therein).

Suppose $x_{i}$, $i=1,\ldots,n$, are i.i.d. continuously distributed
$p$-dimensional random vectors with smooth p.d.f. $f_{0}$ and consider
estimation of the integrated squared density%
\[
\beta_{0}=\int_{\mathbb{R}^{p}}f_{0}(x)^{2}\text{d}x=\mathbb{E}[f_{0}(x_{i})].
\]
A leave-one-out kernel-based estimator is%
\[
\hat{\beta}=\sum_{1\leq i,j\leq n,i\neq j}\mathcal{K}_{h}(x_{i}-x_{j})/n(n-1),
\]
where $\mathcal{K}(u)$ is a symmetric kernel and $\mathcal{K}_{h}%
(u)=h^{-p}\mathcal{K}(u/h)$. This estimator has the V-statistic form of
(\ref{eq:vstat}) with $W_{i}=x_{i}$ and%
\[
\hat{\Gamma}_{n}=1,\text{\hspace{0.2in}}u_{ij}^{n}(W_{i},W_{j})=%
\I
(i\neq j)\{\mathcal{K}_{h}(x_{i}-x_{j})-\beta_{0}\}/\sqrt{n}(n-1).
\]

Let $f_{h}(x)=\int_{\mathbb{R}^{p}}\mathcal{K}(u)f_{0}(x+hu)\text{d}u$ and
$\beta_{h}=\int_{\mathbb{R}^{p}}f_{h}(x)f_{0}(x)$d$x$. By symmetry of
$\mathcal{K}(u)$,%
\[
\mathbb{E}[u_{ij}^{n}(W_{i},W_{j})|W_{i}]=\mathbb{E}[u_{ji}^{n}(W_{j}%
,W_{i})|W_{i}]=\{f_{h}(x_{i})-\beta_{0}\}/\sqrt{n}(n-1),
\]%
\[
\mathbb{E}[u_{ij}^{n}(W_{i},W_{j})]=\{\beta_{h}-\beta_{0}\}/\sqrt{n}(n-1),
\]
so the terms in the decomposition (\ref{eq:decomp}) are of the form%
\[
B_{n}=\sqrt{n}\{\beta_{h}-\beta_{0}\},\text{\qquad}\Psi_{n}=\frac{1}{\sqrt{n}%
}\sum_{1\leq i\leq n}2\{f_{h}(x_{i})-\beta_{h}\},
\]%
\[
U_{n}=\frac{2}{\sqrt{n}(n-1)}\sum_{1\leq i,j\leq n,j<i}\{\mathcal{K}_{h}%
(x_{i}-x_{j})-f_{h}(x_{i})-f_{h}(x_{j})+\beta_{h}\}.
\]

Here, $2\{f_{h}(x_{i})-\beta_{h}\}$ is an approximation to the well known
influence function $2\{f_{0}(x_{i})-\beta_{0}\}$ for estimators of the
integrated squared density. Under regularity conditions, $f_{h}(x_{i})$
converges to $f_{0}(x_{i})$ in mean square as $h\rightarrow0$, so that%
\[
\Psi_{n}=\frac{1}{\sqrt{n}}\sum_{1\leq i\leq n}2\{f_{0}(x_{i})-\beta
_{0}\}+o_{p}(1).
\]
A martingale central limit theorem can be applied as in
\cite{Cattaneo-Crump-Jansson_2014a_ET} to show that the degenerate U-statistic
term $U_{n}$ will be asymptotically normal as $h\rightarrow0$ and
$n\rightarrow\infty$, provided that $n^{2}h^{p}\rightarrow\infty$. It is easy
to show that $n^{2}h^{p}\mathbb{V}[U_{n}]\rightarrow\Delta=\beta_{0}%
\int_{\mathbb{R}^{p}}\mathcal{K}(u)^{2}$d$u$ (under regularity conditions).
Alternative asymptotics occurs when $h^{p}$ shrinks as fast as $1/n$,
resulting in $\mathbb{V}[\Psi_{n}]$ and $\mathbb{V}[U_{n}]$ having the same
magnitude in the limit.

\subsection*{Example 3: \textquotedblleft Many Terms
Asymptotics\textquotedblright}

The previous two examples show how several estimators share the common
structure outlined above. To illustrate how this structure can be applied to
derive new results, the third example studies series estimation in the context
of the partially linear model. The results will shed light on the asymptotic
behavior of this estimator, and the associated inference procedures, when the
number of terms are allowed to grow as fast as the sample size.

Let $(y_{i},x_{i}^{\prime},z_{i}^{\prime})^{\prime}$, $i=1,\ldots,n$, be a
random sample of generated by the partially linear model%
\begin{equation}
y_{i}=x_{i}^{\prime}\beta_{0}+g(z_{i})+\varepsilon_{i},\qquad\mathbb{E}%
[\varepsilon_{i}|x_{i},z_{i}]=0,\label{eq:PLmodel}%
\end{equation}
where $y_{i}$ is a scalar dependent variable, $x_{i}\in\mathbb{R}^{d}$ and
$z_{i}\in\mathbb{R}^{d_{z}}$ are explanatory variables, $\varepsilon_{i}$ is a
disturbance, $g(\cdot)$ is an unknown function, and $\mathbb{E}[\mathbb{V}%
[x_{i}|z_{i}]]$ is of full rank.

A series estimator of $\beta_{0}$ is obtained by regressing $y_{i}$ on $x_{i}$
and approximating functions of $z_{i}$. To describe the estimator, let
$p^{1}(z),$ $p^{2}(z),$ $\ldots$ be approximating functions, such as
polynomials or splines, and let $p_{K}(z)=(p^{1}(z),\ldots,p^{K}(z))^{\prime}$
be a $K$-dimensional vector of such functions. Letting $M_{ij}$ denote the
$(i,j)$-th element of $M=I_{n}-P_{K}(P_{K}^{\prime}P_{K})^{-1}P_{K}^{\prime},$
where $P_{K}=[p_{K}(z_{1}),\ldots,p_{K}(z_{n})]^{\prime}$, a series estimator
of $\beta_{0}$ in (\ref{eq:PLmodel}) is given by%
\[
\hat{\beta}=(\sum_{1\leq i,j\leq n}M_{ij}x_{i}x_{j}^{\prime})^{-1}(\sum_{1\leq
i,j\leq n}M_{ij}x_{i}y_{j}).
\]
\cite{Donald-Newey_1994_JMA} gave conditions for asymptotic normality of this
estimator using standard asymptotics. See also for example \cite{Linton_1995_ECMA}, references therein, for related asymptotic
results when using kernel estimators.

Conditional on $Z=[z_{1},\ldots,z_{n}]^{\prime}$, $\hat{\beta}$ has the
structure outlined earlier:%
\begin{equation}
\sqrt{n}(\hat{\beta}-\beta_{0})=\hat{\Gamma}_{n}^{-1}S_{n}%
,\label{root-n expand}%
\end{equation}
with%
\[
\hat{\Gamma}_{n}=\frac{1}{n}\sum_{1\leq i,j\leq n}M_{ij}x_{i}x_{j}^{\prime
},\qquad S_{n}=\frac{1}{\sqrt{n}}\sum_{1\leq i,j\leq n}x_{i}M_{ij}%
(g_{j}+\varepsilon_{j}),
\]
where $g_{i}=g(z_{i}).$ In other words, $\hat{\beta}$ has the V-statistic form
of (\ref{eq:vstat}) with $W_{i}=(x_{i}^{\prime},\varepsilon_{i})^{\prime}$ and
$u_{ij}^{n}(W_{i},W_{j})=x_{i}M_{ij}(g_{j}+\varepsilon_{j})/\sqrt{n}$.

By $\mathbb{E}[\varepsilon_{i}|x_{i},z_{i}]=0$ we have $\mathbb{E}%
[x_{i}\varepsilon_{i}|Z]=0$. Therefore, letting $u_{ij}^{n}=u_{ij}^{n}%
(W_{i},W_{j})$ as we have done previously, we have
\[
\mathbb{E}[u_{ij}^{n}|Z]=h_{i}M_{ij}g_{j}/\sqrt{n},\text{\qquad}u_{ij}%
^{n}-\mathbb{E}[u_{ij}^{n}|Z]=M_{ij}\left(  v_{i}g_{j}+x_{i}\varepsilon
_{j}\right)  /\sqrt{n},
\]%
\[
\tilde{u}_{ij}^{n}=M_{ij}\left(  v_{j}g_{i}+v_{i}g_{j}+x_{j}\varepsilon
_{i}+x_{i}\varepsilon_{j}\right)  /\sqrt{n},\qquad\mathbb{E}[\tilde{u}%
_{ij}^{n}|W_{i},Z]=M_{ij}\left(  v_{i}g_{j}+h_{j}\varepsilon_{i}\right)
/\sqrt{n},
\]
for $i\neq j$, where $h_{i}=h(z_{i})=\mathbb{E}[x_{i}|z_{i}]$ and $v_{i}%
=x_{i}-h_{i}$. In this case, the bias term in (\ref{eq:decomp}) is%
\[
B_{n}=\frac{1}{\sqrt{n}}\sum_{1\leq i,j\leq n}M_{ij}h_{i}g_{j},
\]
which will be negligible under regularity conditions, as shown in the next
section. Moreover,%
\[
\Psi_{n}=\frac{1}{\sqrt{n}}\sum_{1\leq i\leq n}M_{ii}v_{i}\varepsilon
_{i}+R_{n},\qquad R_{n}=\frac{1}{\sqrt{n}}\sum_{1\leq i,j\leq n}M_{ij}%
(v_{i}g_{j}+h_{i}\varepsilon_{j}),
\]
where $R_{n}$ has mean zero and converges to zero in mean square as $K$ grows,
as further discussed below. Under standard asymptotics $M_{ii}$ will go to one
and hence the limiting variance of the leading term in $\Psi_{n}$ corresponds
to the usual asymptotic variance.

Finally, we find that the degenerate U-statistic term is%
\[
U_{n}=\frac{1}{\sqrt{n}}\sum_{1\leq i,j\leq n,j<i}M_{ij}\left(  v_{i}%
\varepsilon_{j}+v_{j}\varepsilon_{i}\right)  =-\frac{1}{\sqrt{n}}\sum_{1\leq
i,j\leq n,j<i}Q_{ij}\left(  v_{i}\varepsilon_{j}+v_{j}\varepsilon_{i}\right)
.
\]
Remarkably, this term is essentially the same as the degenerate U-statistic
term for JIVE2 that was discussed above. Consequently, the central limit
theorem of \cite{Chao-Swanson-Hausman-Newey-Woutersen_2012_ET} is applicable
to this problem. We will employ it to show that $U_{n}$ is asymptotically
normal as $K\rightarrow\infty$, even when $K/n$ does not converge to zero.

This example highlights a new approach to studying the asymptotic distribution
of semi-linear regression under many terms asymptotics. The alternative
asymptotic approximation is useful, for instance, when the number of
covariates entering the nonparametric part is large relative to the sample
size, as is often the case in empirical applications.

\section{Many Terms Asymptotics\label{section:ManyReg}}

In this section we make precise the discussion given in Example 3, and also
discuss consistent standard error estimation under homoskedasticity. The
estimator $\hat{\beta}$ described in Example 3 can be interpreted as a
two-step semiparametric estimator with tuning parameter $K$, the first step
involving series estimation of the the unknown (regression) functions $g(z)$
and $h(z)$. \cite{Donald-Newey_1994_JMA} gave conditions for asymptotic
normality of this estimator when $K/n\rightarrow0$. Here we generalize their
findings by obtaining an asymptotic distributional result that is valid even
when $K/n$ is bounded away from zero.

The analysis proceeds under the following assumption.

\begin{description}
\item[Assumption PLM] (\textbf{Partially Linear Model})~

\item (a) $(y_{i},x_{i}^{\prime},z_{i}^{\prime})^{\prime}$, $i=1,\ldots,n $,
is a random sample.

\item (b)\ There is a $C<\infty$ such that $\mathbb{E}[\varepsilon_{i}%
^{4}|x_{i},z_{i}]\leq C$ and $\mathbb{E}[\Vert v_{i}\Vert^{4}|z_{i}]\leq C$.

\item (c) There is a $C>0$ such that $\mathbb{E}[\varepsilon_{i}^{2}%
|x_{i},z_{i}]\geq C$ and $\lambda_{\min}(\mathbb{E}[v_{i}v_{i}^{\prime}%
|z_{i}])\geq C$.

\item (d) $\operatorname*{rank}(P_{K})=K$ (a.s.) and there is a $C>0$ such
that $M_{ii}\geq C$.

\item (e) For some $\alpha_{g},\alpha_{h}>0$, there is a $C<\infty$ such that%
\[
\min_{\eta_{g}\in\mathbb{R}^{K}}\mathbb{E}[|g(z_{i})-\eta_{g}^{\prime}%
p_{K}(z_{i})|^{2}]\leq CK^{-2\alpha_{g}},\qquad\min_{\eta_{h}\in
\mathbb{R}^{K\times d}}\mathbb{E}[\Vert h(z_{i})-\eta_{h}^{\prime}p_{K}%
(z_{i})\Vert^{2}]\leq CK^{-2\alpha_{h}}.
\]

\end{description}

Because $\sum_{i=1}^{n}M_{ii}=n-K$, an implication of part (d) is that
$K/n\leq1-C<1$, but crucially Assumption PLM does not imply that
$K/n\rightarrow0$. Part (e) is implied by conventional assumptions from
approximation theory. For instance, when the support of $z_{i}$ is compact
commonly used basis of approximation, such as polynomials or splines, will
satisfy this assumption with $\alpha_{g}=s_{g}/d_{z}$ and $\alpha_{h}%
=s_{h}/d_{z}$, where $s_{g}$ and $s_{h}$ denotes the number of continuous
derivatives of $g(z)$ and $h(z)$, respectively. Further discussion and related
references for several basis of approximation may be found in
\cite{Newey_1997_JoE}, \cite{Chen_2007_Handbook} and
\cite{Belloni-Chernozhukov-Chetverikov-Kato_2015_JoE}, among others.

\subsection{Asymptotic Distribution}

From equation (\ref{root-n expand}), and the discussion in the previous
section, we see that the asymptotic distribution of $\hat{\beta}$ will be
determined by the behavior of $\hat{\Gamma}_{n}$ and $S_{n}$. The following
lemma approximates $\hat{\Gamma}_{n}$ without requiring that $K/n\rightarrow0$.

\begin{lemma}
\label{lemma:GammaHat}If Assumption PLM is satisfied and if $K\rightarrow
\infty$, then%
\[
\hat{\Gamma}_{n}=\Gamma_{n}+o_{p}\left(  1\right)  ,\text{\qquad}\Gamma
_{n}=\frac{1}{n}\sum_{1\leq i\leq n}M_{ii}\mathbb{E}[v_{i}v_{i}^{\prime}%
|z_{i}].
\]

\end{lemma}

Because $\sum_{i=1}^{n}M_{ii}=n-K$, it follows from this result that in the
homoskedastic $v_{i}$ case (i.e., when $\mathbb{E}[v_{i}v_{i}^{\prime}%
|z_{i}]=\mathbb{E}[v_{i}v_{i}^{\prime}]$) $\hat{\Gamma}_{n}$ is close to%
\[
\Gamma_{n}=(1-K/n)\Gamma,\qquad\Gamma=\mathbb{E}[v_{i}v_{i}^{\prime}],
\]
in probability. More generally, with heteroskedasticity, $\hat{\Gamma}_{n}$
will be close to the weighted average $\Gamma_{n}$. Importantly, this result
includes standard asymptotics as a special case when $K/n\rightarrow0 $, where
$\sum_{i=1}^{n}(1-M_{ii})/n=K/n$, the law of large numbers and iterated
expectations imply%
\begin{align*}
\Gamma_{n}  & =\frac{1}{n}\sum_{i=1}^{n}\mathbb{E}[v_{i}v_{i}^{\prime}%
|z_{i}]-\frac{1}{n}\sum_{i=1}^{n}(1-M_{ii})\mathbb{E}[v_{i}v_{i}^{\prime
}|z_{i}]+o_{p}(1)\\
& =\frac{1}{n}\sum_{i=1}^{n}\mathbb{E}[v_{i}v_{i}^{\prime}|z_{i}%
]+o_{p}(1)=\Gamma+o_{p}(1).
\end{align*}

Next, we study
\[
S_{n}=\frac{1}{\sqrt{n}}\sum_{1\leq i,j\leq n}M_{ij}v_{i}\varepsilon_{j}%
+B_{n}+R_{n}.
\]
The following lemma quantifies the magnitude of the bias term $B_{n}$ as well
as the additional variability arising from the (remainder) term $R_{n}$.

\begin{lemma}
\label{lemma:Bias&Remainder}If Assumption PLM is satisfied and if
$K\rightarrow\infty,$ then $B_{n}=O_{p}(\sqrt{n}K^{-\alpha_{g}-\alpha_{h}})$
and $R_{n}=o_{p}(1)$.
\end{lemma}

Like the previous lemma, this lemma does not require $K/n\rightarrow0$.
Interestingly, the bias term $B_{n}$ involves approximation of both unknown
functions $g(z)$ and $h(z)$, implying an implicit trade-off between smoothness
conditions for $g(z)$ and $h(z)$. The implied bias condition $K^{2(\alpha
_{g}+\alpha_{h})}/n\rightarrow\infty$ only requires that $\alpha_{g}%
+\alpha_{h}$ be large enough, but not necessarily that $\alpha_{g}$ and
$\alpha_{h}$ separately be large. It follows that if this bias condition
holds, then%
\[
S_{n}=\frac{1}{\sqrt{n}}\sum_{1\leq i,j\leq n}M_{ij}v_{i}\varepsilon_{j}%
+o_{p}(1),
\]
as claimed in Example 3 above.

Having dispensed with asymptotically negligible contributions to $S_{n}$, we
turn to its leading term. This term is shown below to be asymptotically
Gaussian with asymptotic variance given by%
\[
\Sigma_{n}=\frac{1}{n}\mathbb{V}[\sum_{1\leq i,j\leq n}M_{ij}v_{i}%
\varepsilon_{j}|Z]=\frac{1}{n}\sum_{1\leq i\leq n}M_{ii}^{2}\mathbb{E}%
[v_{i}v_{i}^{\prime}\varepsilon_{i}^{2}|z_{i}]+\frac{1}{n}\sum_{1\leq i,j\leq
n,j\neq i}M_{ij}^{2}\mathbb{E}[v_{i}v_{i}^{\prime}\varepsilon_{j}^{2}%
|z_{i},z_{j}].
\]
Here, the first term following the second equality corresponds to the usual
asymptotic approximation, while the second term adds an additional term that
accounts for large $K$. Once again it is interesting to consider what happens
in some special cases. Under homoskedasticity of $\varepsilon_{i}$ (i.e., when
$\mathbb{E}[\varepsilon_{i}^{2}|x_{i},z_{i}]=\mathbb{E}[\varepsilon_{i}^{2}%
]$),%
\[
\Sigma_{n}=\frac{\sigma_{\varepsilon}^{2}}{n}\sum_{1\leq i,j\leq n}M_{ij}%
^{2}\mathbb{E}[v_{i}v_{i}^{\prime}|z_{i}]=\frac{\sigma_{\varepsilon}^{2}}%
{n}\sum_{1\leq i\leq n}M_{ii}\mathbb{E}[v_{i}v_{i}^{\prime}|z_{i}%
]=\sigma_{\varepsilon}^{2}\Gamma_{n},\qquad\sigma_{\varepsilon}^{2}%
=\mathbb{E}[\varepsilon_{i}^{2}],
\]
because $\sum_{j=1}^{n}M_{ij}^{2}=M_{ii}$. If, in addition, $\mathbb{E}%
[v_{i}v_{i}^{\prime}|z_{i}]=\mathbb{E}[v_{i}v_{i}^{\prime}]$, then $\Sigma
_{n}=\sigma_{\varepsilon}^{2}\left(  1-K/n\right)  \Gamma$. Also, if
$K/n\rightarrow0$, then by $\sum_{1\leq i,j\leq n,i\neq j}M_{ij}^{2}/n\leq
K/n$ and the law of large numbers, we have%
\[
\Sigma_{n}=\frac{1}{n}\sum_{1\leq i\leq n}M_{ii}^{2}\mathbb{E}[v_{i}%
v_{i}^{\prime}\varepsilon_{i}^{2}|z_{i}]+o_{p}\left(  1\right)  =\mathbb{E}%
[v_{i}v_{i}^{\prime}\varepsilon_{i}^{2}]+o_{p}\left(  1\right)  ,
\]
which corresponds to the standard asymptotics limiting variance.

The following theorem combines Lemmas \ref{lemma:GammaHat} and
\ref{lemma:Bias&Remainder} with a central limit theorem for quadratic forms to
show asymptotic normality of $\hat{\beta}$.

\begin{theorem}
\label{thm:AsyNorm}If Assumption PLM is satisfied and if $K^{2(\alpha
_{g}+\alpha_{h})}/n\rightarrow\infty$, then%
\[
\Omega_{n}^{-1/2}\sqrt{n}(\hat{\beta}-\beta_{0})\rightarrow_{d}\mathcal{N}%
(0,I_{d}),\text{\qquad}\Omega_{n}=\Gamma_{n}^{-1}\Sigma_{n}\Gamma_{n}^{-1}.
\]

If, in addition, $\mathbb{E}[\varepsilon_{i}^{2}|x_{i},z_{i}]=\sigma
_{\varepsilon}^{2}$, then $\Omega_{n}=\sigma_{\varepsilon}^{2}\Gamma_{n}^{-1}$.
\end{theorem}

This theorem shows that $\hat{\beta}$ is asymptotically normal when $K/n$ need
not converge to zero. An implication of this result is that inconsistent
series-based nonparametric estimators of the unknown functions $g(z)$ and
$h(z)$ may be employed when forming $\hat{\beta}$, that is, $K/n\nrightarrow0$
is allowed (increasing the variability of the nonparametric estimators),
provided that $K\rightarrow\infty$ (to remove nonparametric smoothing bias).
This asymptotic distributional result does not rely on asymptotic linearity,
nor on the actual convergence of the matrices $\Gamma_{n}$ and $\Sigma_{n}$,
and leads to a new (larger) asymptotic variance that captures terms that are
assumed away by the classical result. The asymptotic distribution result of
\cite{Donald-Newey_1994_JMA} is obtained as a special case where
$K/n\rightarrow0 $. More generally, when $K/n$ does not converge to zero, the
asymptotic variance will be larger than the usual formula because it accounts
for the contribution of \textquotedblleft remainder\textquotedblright\ $U_{n}$
in equation (\ref{eq:decomp}). For instance, when both $\varepsilon_{i}$ and
$v_{i}$ are homoskedastic, the asymptotic variance is%
\[
\Gamma_{n}^{-1}\Sigma_{n}\Gamma_{n}^{-1}=\sigma_{\varepsilon}^{2}\Gamma
_{n}^{-1}=\sigma_{\varepsilon}^{2}\Gamma^{-1}(1-K/n)^{-1}\text{,}%
\]
which is larger than the usual asymptotic variance $\sigma_{\varepsilon}%
^{2}\Gamma^{-1}$ by the degrees of freedom correction $(1-K/n)^{-1}.$

\subsection{Asymptotic Variance Estimation under Homoskedasticity}

Consistent asymptotic variance estimation is useful for large sample
inference. If the assumptions of Theorem \ref{thm:AsyNorm} are satisfied and
if $\hat{\Sigma}_{n}-\Sigma_{n}\rightarrow_{p}0$, then%
\[
\hat{\Omega}_{n}^{-1/2}\sqrt{n}(\hat{\beta}-\beta_{0})\rightarrow
_{d}\mathcal{N}(0,I_{d}),\text{\qquad}\hat{\Omega}_{n}=\hat{\Gamma}_{n}%
^{-1}\hat{\Sigma}_{n}\hat{\Gamma}_{n}^{-1},
\]
implying that valid large-sample confidence intervals and hypothesis tests for
linear and nonlinear transformations of the parameter vector $\beta$ can be
based on $\hat{\Omega}_{n}$.\footnote{Another approach to inference would be
via the bootstrap. For small bandwidth asymptotics,
\cite{Cattaneo-Crump-Jansson_2014b_ET} showed that the standard nonparametric
bootstrap does not provide a valid distributional approximation in general. We
conjecture that the standard nonparametric bootstrap will also fail to provide
valid inference for other alternative asymptotics frameworks.} Under
(conditional)\ heteroskedasticity of unknown form, constructing a consistent
estimator $\hat{\Sigma}_{n}$ turns out to be very challenging if
$K/n\nrightarrow0$. Intuitively, the problem arises because the estimated
residuals entering the construction of $\hat{\Sigma}_{n}$ are not consistent
unless $K/n\rightarrow0$, implying that $\hat{\Sigma}_{n}-\Sigma
_{n}\nrightarrow_{p}0$ in general. Solving this problem is beyond the scope of
this paper. Under homoskedasticity of $\varepsilon_{i}$, however, the
asymptotic variance $\Sigma_{n}$ simplifies and admits a correspondingly
simple consistent estimator. To describe this result, note that if
$\mathbb{E}[\varepsilon_{i}^{2}|x_{i},z_{i}]=\sigma_{\varepsilon}^{2}$ then
$\Sigma_{n}=\sigma_{\varepsilon}^{2}\Gamma_{n}$, where $\hat{\Gamma}%
_{n}-\Gamma_{n}\rightarrow_{p}0$ by Lemma \ref{lemma:GammaHat}. It therefore
suffices to find a consistent estimator of $\sigma_{\varepsilon}^{2}$. Let%
\[
s^{2}=\frac{1}{n-d-K}\sum_{1\leq i\leq n}\hat{\varepsilon}_{i}^{2},\qquad
\hat{\varepsilon}_{i}=\sum_{1\leq j\leq n}M_{ij}(y_{j}-\hat{\beta}^{\prime
}x_{j}),
\]
denote the usual OLS estimator of $\sigma_{\varepsilon}^{2}$ incorporating a
degrees of freedom correction. The following theorem shows that $s^{2}$ is a
consistent estimator, even when the number of terms is \textquotedblleft
large\textquotedblright\ relative to the sample size.

\begin{theorem}
\label{thm:HatAsyVarHOM}Suppose the conditions of Theorem \ref{thm:AsyNorm}
are satisfied. If $\mathbb{E}[\varepsilon_{i}^{2}|x_{i},z_{i}]=\sigma
_{\varepsilon}^{2}$, then $s^{2}\rightarrow_{p}\sigma_{\varepsilon}^{2}$ and
$\hat{\Sigma}_{n}^{\mathtt{HOM}}-\Sigma_{n}\rightarrow_{p}0$, where
$\hat{\Sigma}_{n}^{\mathtt{HOM}}=s^{2}\hat{\Gamma}_{n}$.
\end{theorem}

This theorem provides a distribution free, large sample justification for the
degrees-of-freedom correction required for exact inference under homoskedastic
Gaussian errors. Intuitively, accounting for the correct degrees of freedom is
important whenever the number of terms in the semi-linear model is
\textquotedblleft large\textquotedblright\ relative to the sample size.

\section{Small Simulation Study\label{section:simuls}}

\begin{figure}\centering\includegraphics
[scale=.8] {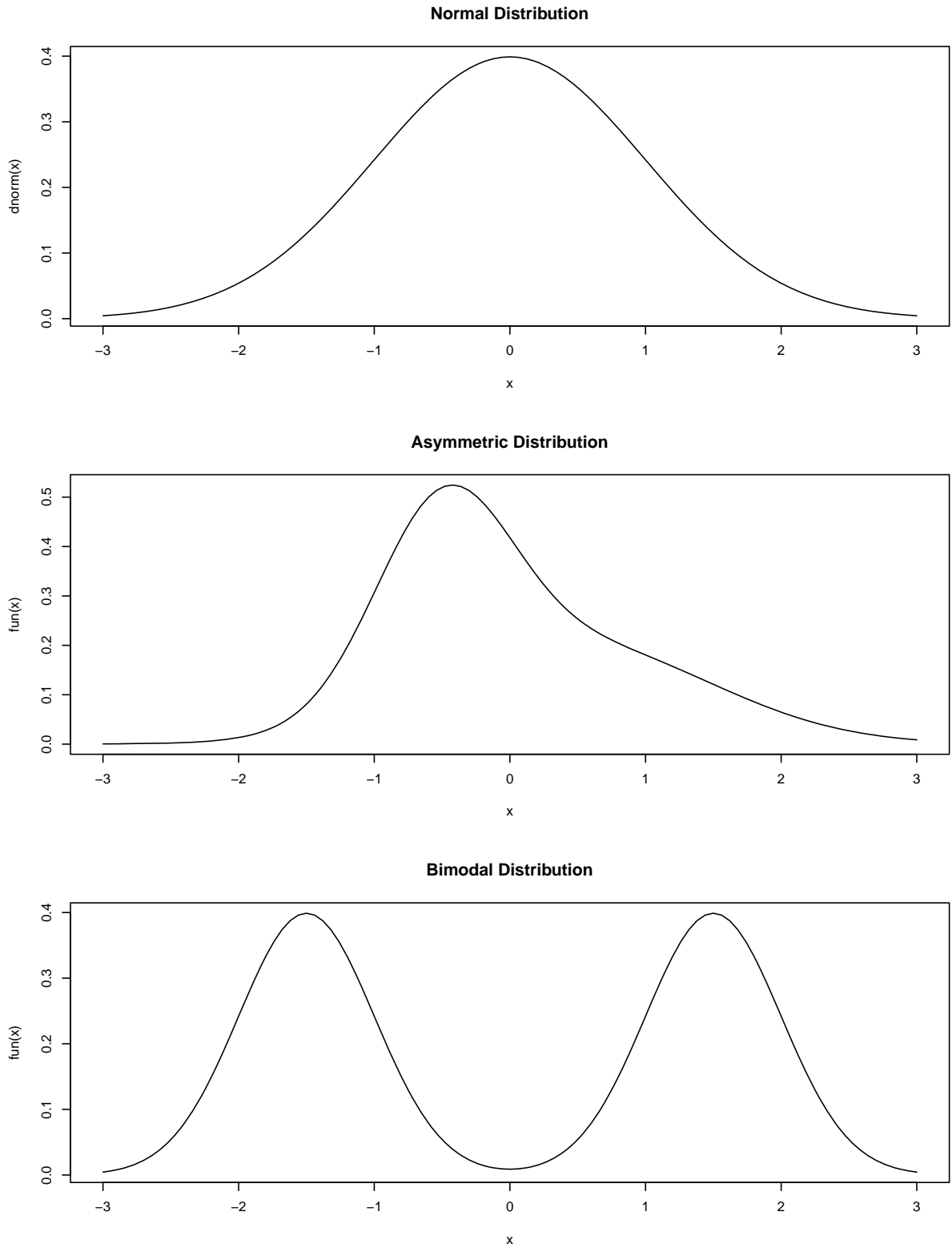}
\caption{Lebesgue Densities of Error Terms Distributions.\label{figure:simuls}%
}
\end{figure}

We conducted a Monte Carlo experiment to explore the extent to which the
asymptotic theoretical results obtained in the previous section are present in
small samples. Using the notation already introduced, we consider the
following partially linear model:%
\[%
\begin{tabular}
[c]{lllll}%
$y_{i}=x_{i}^{\prime}\beta+g(z_{i})+\varepsilon_{i},$ &  & $\mathbb{E}%
[\varepsilon_{i}|x_{i},z_{i}]=0,$ &  & $\mathbb{E}[\varepsilon_{i}^{2}%
|x_{i},z_{i}]=\sigma_{\varepsilon}^{2},$\\
$x_{i}=h(z_{i})+v_{i},$ &  & $\mathbb{E}[v_{i}|z_{i}]=0,$ &  & $\mathbb{E}%
[v_{i}^{2}|z_{i}]=\sigma_{v}^{2}(z_{i}),$%
\end{tabular}
\]
where $d=1$, $\beta=1$, $d_{z}=5$, $z_{i}=(z_{1i},\cdots,z_{d_{z}i})^{\prime}$
with $z_{\ell i}\sim$ i.i.d. $\mathsf{Uniform}(-1,1) $, $\ell=1,\cdots,d_{z}$.
The unknown regression functions are set to $g(z_{i})=h(z_{i})=\exp(\Vert
z_{i}\Vert^{2})$, which are not additive separable in the covariates $z_{i}$.
The simulation study is based on $S=5,000$ replications, each replication
taking a random sample of size $n=500 $ with all random variables generated
independently. We consider $6$ data generating processes (DGPs) as follows:%
\[%
\begin{tabular}
[c]{cccc}%
\multicolumn{4}{c}{Data Generating Process for Monte Carlo Experiment}%
\\\hline\hline
& \multicolumn{3}{c}{$(\varepsilon_{i},v_{i})$ -- Distributions}\\\cline{2-4}
& Gaussian & Asymmetric & Bimodal\\\hline
\multicolumn{1}{l}{$\sigma_{v}^{2}(z_{i})=1$} & Model 1 & Model 3 & Model 5\\
\multicolumn{1}{l}{$\sigma_{v}^{2}(z_{i})=\varsigma(1+\Vert z_{i}\Vert
^{2})^{2}$} & Model 2 & Model 4 & Model 6\\\hline
\end{tabular}
\]
Specifically, Models 1, 3 and 5 correspond to homoskedastic (in $v_{i}$) DGPs,
while Models 2, 4 and 5 correspond to heteroskedastic (in $v_{i}$) DGPs. For
the latter models, the constant $\varsigma$ was chosen so that $\mathbb{E}%
[v_{i}^{2}]=1$. The three distributions considered for the unobserved error
terms $\varepsilon_{i}$ and $v_{i}$ are: the standard Normal (labelled
\textquotedblleft Gaussian\textquotedblright) and two Mixture of Normals
inducing either an asymmetric or a bimodal distribution; their Lebesgue
densities are depicted in Figure
\ref{figure:simuls}%
. We explored other specifications for the regression functions,
heteroskedasticity form, and distributional assumptions, but we do not report
these additional results because they were qualitative similar to those
discussed here.

The estimators considered in the Monte Carlo experiment are constructed using
power series approximations. We do not impose additive separability on the
basis, though we do restrict the interaction terms to not exceed degree 5. To
be specific, we consider the following polynomial basis expansion:%
\[%
\begin{tabular}
[c]{ccc}%
\multicolumn{3}{c}{Polynomial Basis Expansion: $d_{z}=5$ and $n=500$%
}\\\hline\hline
$K$ & $p_{K}(z_{i})$ & $K/n$\\\hline
$6$ & $(1,z_{1i},z_{2i},z_{3i},z_{4i},z_{5i})^{\prime}$ & $0.012$\\
$11$ & $(p_{6}(z_{i})^{\prime},z_{1i}^{2},z_{2i}^{2},z_{3i}^{2},z_{4i}%
^{2},z_{5i}^{2})^{\prime}$ & $0.022 $\\
$21$ & $p_{11}(z_{i})$ $+$ first-order interactions & $0.042$\\
$26$ & $(p_{21}(z_{i})^{\prime},z_{1i}^{3},z_{2i}^{3},z_{3i}^{3},z_{4i}%
^{3},z_{5i}^{3})^{\prime}$ & $0.052 $\\
$56$ & $p_{26}(z_{i})$ $+$ second-order interactions & $0.112$\\
$61$ & $(p_{56}(z_{i})^{\prime},z_{1i}^{4},z_{2i}^{4},z_{3i}^{4},z_{4i}%
^{4},z_{5i}^{4})^{\prime}$ & $0.122 $\\
$126$ & $p_{61}(z_{i})$ $+$ third-order interactions & $0.252$\\
$131$ & $(p_{126}(z_{i})^{\prime},z_{1i}^{5},z_{2i}^{5},z_{3i}^{5},z_{4i}%
^{5},z_{5i}^{5})^{\prime}$ & $0.262 $\\
$252$ & $p_{131}(z_{i})$ $+$ fourth-order interactions & $0.504$\\
$257$ & $(p_{252}(z_{i})^{\prime},z_{1i}^{6},z_{2i}^{6},z_{3i}^{6},z_{4i}%
^{6},z_{5i}^{6})^{\prime}$ & $0.514 $\\
$262$ & $(p_{257}(z_{i})^{\prime},z_{1i}^{7},z_{2i}^{7},z_{3i}^{7},z_{4i}%
^{7},z_{5i}^{7})^{\prime}$ & $0.524 $\\
$267$ & $(p_{262}(z_{i})^{\prime},z_{1i}^{8},z_{2i}^{8},z_{3i}^{8},z_{4i}%
^{8},z_{5i}^{8})^{\prime}$ & $0.534 $\\
$272$ & $(p_{267}(z_{i})^{\prime},z_{1i}^{9},z_{2i}^{9},z_{3i}^{9},z_{4i}%
^{9},z_{5i}^{9})^{\prime}$ & $0.544 $\\
$277$ & $(p_{272}(z_{i})^{\prime},z_{1i}^{10},z_{2i}^{10},z_{3i}^{10}%
,z_{4i}^{10},z_{5i}^{10})^{\prime}$ & $0.554$\\\hline
\end{tabular}
\]

Thus, our simulations explore the consequences of introducing many terms in
the partially linear model by varying $K$ on the grid above from $K=6$ to
$K=277$, which gives a range for $K/n$ of $\{0.012,\cdots,0.554\}$. For each
point on the grid of $K/n$, we report average bias, average standard
deviation, mean square error and average standarized bias of $\hat{\beta}$
across simulations. We also consider the coverage error rates and interval
length for two asymptotic $95\%$ confidence intervals:%
\[
\text{CI}_{0}=\left[  \hat{\beta}-\Phi_{1-\alpha/2}^{-1}\frac{\hat{\sigma}%
\hat{\Gamma}_{n}^{-1/2}}{\sqrt{n}}\quad,\quad\hat{\beta}+\Phi_{1-\alpha
/2}^{-1}\frac{\hat{\sigma}\hat{\Gamma}_{n}^{-1/2}}{\sqrt{n}}\right]  ,
\]%
\[
\text{CI}_{1}=\left[  \hat{\beta}-\Phi_{1-\alpha/2}^{-1}\frac{s\hat{\Gamma
}_{n}^{-1/2}}{\sqrt{n}}\quad,\quad\hat{\beta}+\Phi_{1-\alpha/2}^{-1}%
\frac{s\hat{\Gamma}_{n}^{-1/2}}{\sqrt{n}}\right]  ,
\]
where $\hat{\sigma}^{2}=(n-d-K)s^{2}/n$, and $\Phi_{u}^{-1}=\Phi^{-1}(u)$
denotes the inverse of the Gaussian distribution function. That is, CI$_{0}$
and CI$_{1}$ are formed employing the t-statistic constructed using the
homoskedasticity-consistent variance estimators without and with degrees of
freedom correction, respectively.

The main findings from the Monte Carlo experiment are presented in Tables
\ref{table:simuls1}%
--%
\ref{table:simuls3}%
. All results are consistent with the theoretical conclusions presented in the
previous section. First, the results for standard Normal and non-Normal errors
are qualitatively similar. This indicates that the Gaussian approximation
obtained in Theorem \ref{thm:AsyNorm} is a good approximation in finite
samples, even when $K$ is a nontrivial fraction of the sample size. Second, as
expected, a small choice of $K$ leads to important smoothing biases. This
affects the finite sample properties of the point estimators as well as the
distributional approximations obtained in this paper. In particular, it
affects the empirical size of all the confidence intervals. Third, in all
cases the results under homoskedasticity or heteroskedasticity in $v_{i}$ are
qualitatively similar, showing that our theoretical results provide a good
finite sample approximation in both cases, even when $K$ is a nontrivial
fraction of the sample size. Fourth, as suggested by Theorem
\ref{thm:HatAsyVarHOM}, confidence intervals without degrees of freedom
correction (CI$_{0}$) are under-sized, while the analogue confidence intervals
with degrees of freedom correction (CI$_{1}$) have close-to-correct empirical
size in all cases. This result shows that the degrees of freedom correction is
crucial to achieve close-to-correct empirical size when $K/n$ is non-negligible.

In conclusion, we found in our small-scale simulation study that our
theoretical results for the partially linear model with possibly many terms
provide good approximation in samples of moderate size. In particular, under
homoskedasticity of $\varepsilon_{i}$, we showed that confidence intervals
constructed using $s^{2}$ exhibit good empirical coverage even when $K/n$ is
\textquotedblleft large\textquotedblright. We also confirmed that the Gaussian
distributional approximation given in Theorem \ref{thm:AsyNorm} represents
well the finite sample distribution of $\hat{\beta}$ even when $K/n$ is
\textquotedblleft large\textquotedblright.

\section{Conclusion\label{section:conclusion}}

This paper showed that the many instrument asymptotics and the small bandwidth
asymptotics shared a common structure based on a V-statistic, with a remainder
term that is asymptotically normal when the number of term diverges to
infinity or the bandwidth shrinks to zero. This feature is particularly useful
to obtain new results for other semiparametric estimators. In this paper we
employ this common structure to derive a new alternative large-sample
distributional approximation for a series estimator of the partially linear
model, which implied a new (larger) asymptotic variance formula.

Our results apply to a class of semiparametric estimators $\hat{\beta}$
satisfying%
\[
\sqrt{n}(\hat{\beta}-\beta_{0})=\hat{\Gamma}_{n}^{-1}S_{n}+o_{p}(1),
\]
where $\hat{\Gamma}_{n}$ and $S_{n}$ take a particular V-stastistic form, as
discussed in Section \ref{section:CommonStructure}. This class of
semiparametric estimators covers several interesting problems, but it is by no
means exhaustive. For example, \cite{Cattaneo-Jansson_2015_wp-Kernels} show
that a large class of (kernel-based) semiparametric estimators admit an
expansion of the form%
\[
\sqrt{n}(\hat{\beta}-\beta_{0})=\hat{\Gamma}_{n}^{-1}S_{n}-\mathcal{B}%
_{n}+o_{p}(1),
\]
where the bias term $\mathcal{B}_{n}$ is quantitatively and conceptually
distinct from the smoothing bias $B_{n}$ described in Section
\ref{section:CommonStructure} and, crucially, dominates the quadratic term
$U_{n} $ arising from the V-statistic $S_{n}$; that is, $U_{n}=o_{p}%
(\mathcal{B}_{n})$ in that setting. Nevertheless, the structure we have
considered in this paper is useful, providing new results for the partially
linear model and a common structure for disparate literatures on many
instruments and small bandwidths.

\section{Appendix: Proofs}

All statements involving conditional expectations are understood to hold
almost surely. Qualifiers such as \textquotedblleft a.s.\textquotedblright%
\ will be omitted to conserve space. Throughout the appendix, $C$ will denote
a generic constant that may take different values in each case.\bigskip

\textbf{Proof of Lemma \ref{lemma:GammaHat}}. Let $X=[x_{1},\ldots
,x_{n}]^{\prime}$, $H=[h_{1},\ldots,h_{n}]^{\prime}$, and $V=[v_{1}%
,\ldots,v_{n}]^{\prime}$. By Assumption PLM and the Markov inequality,%
\[
\operatorname*{tr}(\frac{1}{n}H^{\prime}MH)=\min_{\eta_{h}\in\mathbb{R}%
^{K\times d}}\frac{1}{n}\sum_{1\leq i\leq n}\Vert h(z_{i})-\eta_{h}^{\prime
}p_{K}(z_{i})\Vert^{2}=O_{p}(K^{-2\alpha_{h}})\rightarrow_{p}0.
\]
Also, $V^{\prime}V/n=O_{p}(1)$ by Assumption PLM and the Markov inequality, so
by the Cauchy-Schwarz inequality and $M$ idempotent, $\Vert H^{\prime
}MV/n\Vert\leq\lbrack\operatorname*{tr}(H^{\prime}MH/n)\operatorname*{tr}%
(V^{\prime}V/n)]^{1/2}\rightarrow_{p}0.$ By the triangle inequality, we then
have%
\[
\hat{\Gamma}_{n}=\frac{1}{n}X^{\prime}MX=\frac{1}{n}(V+H)^{\prime}%
M(V+H)=\frac{1}{n}V^{\prime}MV+o_{p}(1).
\]
Next, by Lemma A1 of \cite{Chao-Swanson-Hausman-Newey-Woutersen_2012_ET},%
\[
\frac{1}{n}V^{\prime}MV=\frac{1}{n}\sum_{1\leq i\leq n}M_{ii}v_{i}%
v_{i}^{\prime}+\frac{1}{n}\sum_{1\leq i,j\leq n,j\neq i}M_{ij}v_{i}%
v_{j}^{\prime}=\frac{1}{n}\sum_{1\leq i\leq n}M_{ii}v_{i}v_{i}^{\prime}%
+o_{p}(1).
\]
Finally, by the Markov inequality and using $\mathbb{E}[n^{-1}\sum_{1\leq
i\leq n}M_{ii}v_{i}v_{i}^{\prime}|Z]=\Gamma_{n}$,%
\[
\frac{1}{n}\sum_{1\leq i\leq n}M_{ii}v_{i}v_{i}^{\prime}-\Gamma_{n}%
\rightarrow_{p}0
\]
because Assumption PLM implies that $v_{i}v_{i}^{\prime}$ and $v_{j}%
v_{j}^{\prime}$ are uncorrelated conditional on $Z$ and that $\mathbb{E}%
[M_{ii}^{2}\Vert v_{i}\Vert^{4}|Z]\leq C$.\qquad$\blacksquare$\bigskip

\textbf{Proof of Lemma \ref{lemma:Bias&Remainder}}. Let $G=[g_{1},\ldots
,g_{n}]^{\prime}$ and $\varepsilon=[\varepsilon_{1},\ldots,\varepsilon
_{n}]^{\prime}$. By the Cauchy-Schwarz inequality, $M$ idempotent, Assumption
PLM, and the Markov inequality,%

\[
\Vert\frac{1}{n}G^{\prime}MH\Vert\leq\sqrt{\operatorname*{tr}(\frac{1}%
{n}G^{\prime}MG)}\sqrt{\operatorname*{tr}(\frac{1}{n}H^{\prime}MH)}%
=O_{p}(K^{-\alpha_{g}-\alpha_{h}}),
\]
which gives $B_{n}=G^{\prime}MH/\sqrt{n}=O_{p}(\sqrt{n}K^{-\alpha_{g}%
-\alpha_{h}})$.

Also, $R_{n}=(V^{\prime}MG+H^{\prime}M\varepsilon)/\sqrt{n}=O_{p}%
(K^{-\alpha_{g}}+K^{-\alpha_{h}})=o_{p}(1)$ because%
\[
\mathbb{E}[\Vert\frac{1}{\sqrt{n}}V^{\prime}MG\Vert^{2}|Z]=\frac{1}%
{n}G^{\prime}M\mathbb{E}[VV^{\prime}|Z]MG\leq C\frac{1}{n}G^{\prime}%
MG=O_{p}(K^{-2\alpha_{g}})
\]
and%
\[
\mathbb{E}[\Vert\frac{1}{\sqrt{n}}H^{\prime}M\varepsilon\Vert^{2}%
|Z]=\operatorname*{tr}(\frac{1}{n}H^{\prime}M\mathbb{E}[\varepsilon
\varepsilon^{\prime}|Z]MH)\leq C\operatorname*{tr}(\frac{1}{n}H^{\prime
}MH)=O_{p}(K^{-2\alpha_{h}})
\]
by Assumption PLM and the Markov inequality.\qquad$\blacksquare$\bigskip

\textbf{Proof of Theorem \ref{thm:AsyNorm}}. By Lemma A2 of
\cite{Chao-Swanson-Hausman-Newey-Woutersen_2012_ET},
\[
\Sigma_{n}^{-1/2}\frac{1}{\sqrt{n}}\sum_{1\leq i,j\leq n}M_{ij}v_{i}%
\varepsilon_{j}\rightarrow_{d}\mathcal{N}(0,I_{d})
\]
under Assumption PLM. Combining this result with Lemmas \ref{lemma:GammaHat}
and \ref{lemma:Bias&Remainder}, we obtain the results stated in the
theorem.\qquad$\blacksquare$\medskip

\textbf{Proof of Theorem \ref{thm:HatAsyVarHOM}}. Let $Y=[y_{1},\ldots,y_{n}]$
and $\hat{\varepsilon}=[\hat{\varepsilon}_{1},\ldots,\hat{\varepsilon}%
_{n}]^{\prime}=M(Y-X\hat{\beta})$. It follows similarly to the proof of Lemma
\ref{lemma:GammaHat} that%
\begin{align*}
\frac{1}{n}\varepsilon^{\prime}M\varepsilon & =\frac{1}{n}\sum_{1\leq i\leq
n}M_{ii}\varepsilon_{i}^{2}+\frac{1}{n}\sum_{1\leq i,j\leq n,j\neq
i}\varepsilon_{i}M_{ij}\varepsilon_{j}\\
& =\frac{1}{n}\sum_{1\leq i\leq n}M_{ii}\mathbb{E}[\varepsilon_{i}^{2}%
|z_{i}]+o_{p}\left(  1\right)  =\frac{n-K}{n}\sigma_{\varepsilon}^{2}%
+o_{p}(1),
\end{align*}
so it suffices to show that $\hat{\varepsilon}^{\prime}\hat{\varepsilon
}/n=\varepsilon^{\prime}M\varepsilon/n+o_{p}(1)$.

Lemma \ref{lemma:GammaHat} and $\hat{\beta}-\beta=o_{p}(1)$ imply $(\hat
{\beta}-\beta)^{\prime}X^{\prime}MX(\hat{\beta}-\beta)/n=o_{p}\left(
1\right)  $, which together with the Cauchy-Schwarz inequality and
$\varepsilon^{\prime}M\varepsilon/n=O_{p}(1)$ gives%
\begin{align*}
\frac{1}{n}(Y-X\hat{\beta}-G)^{\prime}M(Y-X\hat{\beta}-G) &  =\frac{1}%
{n}\varepsilon^{\prime}M\varepsilon+\frac{1}{n}(\hat{\beta}-\beta)^{\prime
}X^{\prime}MX(\hat{\beta}-\beta)-\frac{1}{n}2\varepsilon^{\prime}MX(\hat
{\beta}-\beta)\\
&  =\frac{1}{n}\varepsilon^{\prime}M\varepsilon+o_{p}(1).
\end{align*}
Similarly, $G^{\prime}MG/n=o_{p}\left(  1\right)  $ together with
$(Y-X\hat{\beta}-G)^{\prime}M(Y-X\hat{\beta}-G)/n=O_{p}\left(  1\right)  $ and
the Cauchy-Schwarz inequality gives%
\[
\frac{1}{n}\hat{\varepsilon}^{\prime}\hat{\varepsilon}=\frac{1}{n}%
(Y-X\hat{\beta})^{\prime}M(Y-X\hat{\beta})=\frac{1}{n}(Y-X\hat{\beta
}-G)^{\prime}M(Y-X\hat{\beta}-G)+o_{p}(1).
\]
The conclusion follows by the triangle inequality.\qquad$\blacksquare$

\bibliographystyle{econometrica}
\bibliography{altasym_19Dec2014}

\newpage%

\newpage%

\begin{table}\centering\vspace{-.2in}
\caption{Simulation Results, Models $1-2$, Gaussian Distribution.\label
{table:simuls1}}
\subfloat[Model 1: Homoskedastic $v_i$]{
\begin{tabular}{rrrrrrrrr}
\hline\hline
\multicolumn{1}{c}{$K/n$}&\multicolumn{1}{c}{$\mathsf{Bias}$}&\multicolumn{1}{c}{$\mathsf{SD}$}&\multicolumn{1}{c}{$\mathsf{RMSE}$}&\multicolumn{1}{c}{$\frac{\mathsf{Bias}}{\mathsf{SD}}$}&\multicolumn{1}{c}{CI$_0$}&\multicolumn{1}{c}{CI$_1$}&\multicolumn{1}{c}{$\hat\sigma$}&\multicolumn{1}{c}{$s$}\tabularnewline
\hline
$0.012$&$0.481$&$0.040$&$0.483$&$11.898$&$0.000$&$0.000$&$0.039$&$0.039$\tabularnewline
$0.022$&$0.001$&$0.045$&$0.045$&$ 0.031$&$0.947$&$0.950$&$0.045$&$0.045$\tabularnewline
$0.042$&$0.002$&$0.047$&$0.047$&$ 0.051$&$0.939$&$0.945$&$0.045$&$0.046$\tabularnewline
$0.052$&$0.002$&$0.046$&$0.046$&$ 0.049$&$0.940$&$0.947$&$0.045$&$0.046$\tabularnewline
$0.112$&$0.002$&$0.047$&$0.047$&$ 0.041$&$0.936$&$0.952$&$0.045$&$0.048$\tabularnewline
$0.122$&$0.000$&$0.048$&$0.048$&$ 0.005$&$0.935$&$0.949$&$0.045$&$0.048$\tabularnewline
$0.252$&$0.001$&$0.052$&$0.052$&$ 0.013$&$0.907$&$0.947$&$0.045$&$0.052$\tabularnewline
$0.262$&$0.000$&$0.052$&$0.052$&$-0.008$&$0.904$&$0.949$&$0.045$&$0.052$\tabularnewline
$0.504$&$0.000$&$0.063$&$0.063$&$ 0.003$&$0.841$&$0.951$&$0.045$&$0.064$\tabularnewline
$0.514$&$0.000$&$0.064$&$0.064$&$-0.002$&$0.828$&$0.947$&$0.045$&$0.064$\tabularnewline
$0.524$&$0.000$&$0.064$&$0.064$&$-0.003$&$0.827$&$0.948$&$0.045$&$0.065$\tabularnewline
$0.534$&$0.000$&$0.066$&$0.066$&$-0.003$&$0.821$&$0.950$&$0.045$&$0.066$\tabularnewline
$0.544$&$0.001$&$0.068$&$0.068$&$ 0.010$&$0.803$&$0.946$&$0.045$&$0.067$\tabularnewline
$0.554$&$0.000$&$0.067$&$0.067$&$ 0.004$&$0.808$&$0.949$&$0.045$&$0.067$\tabularnewline
\hline
\end{tabular}

}\qquad\subfloat
[Model 2: Heteroskedastic $v_i$]{
\begin{tabular}{rrrrrrrrr}
\hline\hline
\multicolumn{1}{c}{$K/n$}&\multicolumn{1}{c}{$\mathsf{Bias}$}&\multicolumn{1}{c}{$\mathsf{SD}$}&\multicolumn{1}{c}{$\mathsf{RMSE}$}&\multicolumn{1}{c}{$\frac{\mathsf{Bias}}{\mathsf{SD}}$}&\multicolumn{1}{c}{CI$_0$}&\multicolumn{1}{c}{CI$_1$}&\multicolumn{1}{c}{$\hat\sigma$}&\multicolumn{1}{c}{$s$}\tabularnewline
\hline
$0.012$&$0.483$&$0.046$&$0.485$&$10.460$&$0.000$&$0.000$&$0.039$&$0.040$\tabularnewline
$0.022$&$0.002$&$0.045$&$0.045$&$ 0.034$&$0.949$&$0.953$&$0.045$&$0.046$\tabularnewline
$0.042$&$0.001$&$0.046$&$0.046$&$ 0.015$&$0.946$&$0.949$&$0.045$&$0.046$\tabularnewline
$0.052$&$0.002$&$0.046$&$0.046$&$ 0.034$&$0.947$&$0.955$&$0.045$&$0.046$\tabularnewline
$0.112$&$0.001$&$0.049$&$0.049$&$ 0.015$&$0.932$&$0.950$&$0.045$&$0.048$\tabularnewline
$0.122$&$0.001$&$0.049$&$0.049$&$ 0.025$&$0.929$&$0.946$&$0.045$&$0.049$\tabularnewline
$0.252$&$0.000$&$0.052$&$0.052$&$ 0.009$&$0.914$&$0.951$&$0.046$&$0.053$\tabularnewline
$0.262$&$0.001$&$0.053$&$0.053$&$ 0.025$&$0.915$&$0.952$&$0.046$&$0.054$\tabularnewline
$0.504$&$0.000$&$0.068$&$0.068$&$ 0.002$&$0.827$&$0.947$&$0.048$&$0.068$\tabularnewline
$0.514$&$0.001$&$0.068$&$0.068$&$ 0.019$&$0.829$&$0.953$&$0.048$&$0.068$\tabularnewline
$0.524$&$0.003$&$0.068$&$0.069$&$ 0.050$&$0.824$&$0.953$&$0.047$&$0.069$\tabularnewline
$0.534$&$0.000$&$0.070$&$0.070$&$ 0.003$&$0.819$&$0.949$&$0.048$&$0.070$\tabularnewline
$0.544$&$0.002$&$0.070$&$0.070$&$ 0.024$&$0.819$&$0.948$&$0.048$&$0.071$\tabularnewline
$0.554$&$0.000$&$0.074$&$0.074$&$-0.004$&$0.801$&$0.943$&$0.048$&$0.072$\tabularnewline
\hline
\end{tabular}

}\\
\footnotesize\raggedright Notes:\newline(i) columns $\mathsf{Bias}%
$, $\mathsf{SD}$, $\mathsf{RMSE}$ and $\frac{\mathsf{Bias}}{\mathsf{SD}%
}%
$ report, respectively, average bias, average standard deviation, root mean square error, and average standarized bias of the estimator $\hat
{\beta}$ across simulations;\newline
(ii) columns CI$_0$ and CI$_1$ report empirical coverage for homoskedastic-consistent confidence intervals, respectively, without and with degrees of freedom correction;\newline
(iii) columns $\hat\sigma
$ and $s$ report the average across simulations of the standard errors estimators, respectively, without and with degrees of freedom correction.
\end{table}%
%

\begin{table}\centering\vspace{-.2in}
\caption{Simulation Results, Models $3-4$, Asymmetric Distribution.\label
{table:simuls2}}
\subfloat[Model 3: Homoskedastic $v_i$]{
\begin{tabular}{rrrrrrrrr}
\hline\hline
\multicolumn{1}{c}{$K/n$}&\multicolumn{1}{c}{$\mathsf{Bias}$}&\multicolumn{1}{c}{$\mathsf{SD}$}&\multicolumn{1}{c}{$\mathsf{RMSE}$}&\multicolumn{1}{c}{$\frac{\mathsf{Bias}}{\mathsf{SD}}$}&\multicolumn{1}{c}{CI$_0$}&\multicolumn{1}{c}{CI$_1$}&\multicolumn{1}{c}{$\hat\sigma$}&\multicolumn{1}{c}{$s$}\tabularnewline
\hline
$0.012$&$0.481$&$0.039$&$0.483$&$12.486$&$0.000$&$0.000$&$0.038$&$0.038$\tabularnewline
$0.022$&$0.002$&$0.043$&$0.043$&$ 0.040$&$0.943$&$0.946$&$0.042$&$0.042$\tabularnewline
$0.042$&$0.001$&$0.044$&$0.044$&$ 0.032$&$0.942$&$0.947$&$0.042$&$0.043$\tabularnewline
$0.052$&$0.001$&$0.043$&$0.043$&$ 0.023$&$0.946$&$0.954$&$0.042$&$0.043$\tabularnewline
$0.112$&$0.001$&$0.045$&$0.045$&$ 0.023$&$0.931$&$0.947$&$0.042$&$0.044$\tabularnewline
$0.122$&$0.002$&$0.045$&$0.045$&$ 0.036$&$0.936$&$0.951$&$0.042$&$0.045$\tabularnewline
$0.252$&$0.001$&$0.049$&$0.049$&$ 0.013$&$0.902$&$0.950$&$0.042$&$0.048$\tabularnewline
$0.262$&$0.001$&$0.049$&$0.049$&$ 0.013$&$0.915$&$0.953$&$0.042$&$0.049$\tabularnewline
$0.504$&$0.000$&$0.060$&$0.060$&$ 0.001$&$0.829$&$0.950$&$0.042$&$0.059$\tabularnewline
$0.514$&$0.000$&$0.060$&$0.060$&$-0.007$&$0.828$&$0.948$&$0.042$&$0.060$\tabularnewline
$0.524$&$0.000$&$0.060$&$0.060$&$-0.006$&$0.830$&$0.952$&$0.042$&$0.061$\tabularnewline
$0.534$&$0.000$&$0.061$&$0.061$&$-0.001$&$0.819$&$0.950$&$0.042$&$0.061$\tabularnewline
$0.544$&$0.000$&$0.062$&$0.062$&$ 0.000$&$0.809$&$0.951$&$0.042$&$0.062$\tabularnewline
$0.554$&$0.001$&$0.064$&$0.064$&$ 0.009$&$0.794$&$0.944$&$0.042$&$0.063$\tabularnewline
\hline
\end{tabular}

}\qquad\subfloat
[Model 4: Heteroskedastic $v_i$]{
\begin{tabular}{rrrrrrrrr}
\hline\hline
\multicolumn{1}{c}{$K/n$}&\multicolumn{1}{c}{$\mathsf{Bias}$}&\multicolumn{1}{c}{$\mathsf{SD}$}&\multicolumn{1}{c}{$\mathsf{RMSE}$}&\multicolumn{1}{c}{$\frac{\mathsf{Bias}}{\mathsf{SD}}$}&\multicolumn{1}{c}{CI$_0$}&\multicolumn{1}{c}{CI$_1$}&\multicolumn{1}{c}{$\hat\sigma$}&\multicolumn{1}{c}{$s$}\tabularnewline
\hline
$0.012$&$ 0.485$&$0.046$&$0.488$&$10.566$&$0.000$&$0.000$&$0.038$&$0.038$\tabularnewline
$0.022$&$ 0.001$&$0.042$&$0.042$&$ 0.031$&$0.947$&$0.949$&$0.042$&$0.043$\tabularnewline
$0.042$&$ 0.001$&$0.043$&$0.043$&$ 0.025$&$0.946$&$0.951$&$0.042$&$0.043$\tabularnewline
$0.052$&$ 0.002$&$0.044$&$0.044$&$ 0.047$&$0.937$&$0.943$&$0.042$&$0.043$\tabularnewline
$0.112$&$ 0.002$&$0.045$&$0.045$&$ 0.037$&$0.933$&$0.945$&$0.043$&$0.045$\tabularnewline
$0.122$&$ 0.001$&$0.046$&$0.046$&$ 0.025$&$0.929$&$0.945$&$0.043$&$0.046$\tabularnewline
$0.252$&$ 0.000$&$0.050$&$0.050$&$-0.004$&$0.910$&$0.949$&$0.043$&$0.050$\tabularnewline
$0.262$&$ 0.001$&$0.050$&$0.050$&$ 0.020$&$0.907$&$0.951$&$0.043$&$0.050$\tabularnewline
$0.504$&$ 0.000$&$0.064$&$0.064$&$-0.002$&$0.832$&$0.947$&$0.045$&$0.064$\tabularnewline
$0.514$&$ 0.001$&$0.065$&$0.065$&$ 0.008$&$0.827$&$0.948$&$0.045$&$0.064$\tabularnewline
$0.524$&$-0.001$&$0.065$&$0.065$&$-0.015$&$0.817$&$0.948$&$0.045$&$0.065$\tabularnewline
$0.534$&$ 0.001$&$0.066$&$0.066$&$ 0.013$&$0.824$&$0.948$&$0.045$&$0.065$\tabularnewline
$0.544$&$ 0.000$&$0.067$&$0.067$&$-0.002$&$0.799$&$0.951$&$0.045$&$0.066$\tabularnewline
$0.554$&$ 0.000$&$0.067$&$0.067$&$-0.001$&$0.811$&$0.948$&$0.045$&$0.067$\tabularnewline
\hline
\end{tabular}

}\\
\footnotesize\raggedright Notes:\newline(i) columns $\mathsf{Bias}%
$, $\mathsf{SD}$, $\mathsf{RMSE}$ and $\frac{\mathsf{Bias}}{\mathsf{SD}%
}%
$ report, respectively, average bias, average standard deviation, root mean square error, and average standarized bias of the estimator $\hat
{\beta}$ across simulations;\newline
(ii) columns CI$_0$ and CI$_1$ report empirical coverage for homoskedastic-consistent confidence intervals, respectively, without and with degrees of freedom correction;\newline
(iii) columns $\hat\sigma
$ and $s$ report the average across simulations of the standard errors estimators, respectively, without and with degrees of freedom correction.
\end{table}%
%

\begin{table}\centering\vspace{-.2in}
\caption{Simulation Results, Models $5-6$, Bimodal Distribution.\label
{table:simuls3}}
\subfloat[Model 5: Homoskedastic $v_i$]{
\begin{tabular}{rrrrrrrrr}
\hline\hline
\multicolumn{1}{c}{$K/n$}&\multicolumn{1}{c}{$\mathsf{Bias}$}&\multicolumn{1}{c}{$\mathsf{SD}$}&\multicolumn{1}{c}{$\mathsf{RMSE}$}&\multicolumn{1}{c}{$\frac{\mathsf{Bias}}{\mathsf{SD}}$}&\multicolumn{1}{c}{CI$_0$}&\multicolumn{1}{c}{CI$_1$}&\multicolumn{1}{c}{$\hat\sigma$}&\multicolumn{1}{c}{$s$}\tabularnewline
\hline
$0.012$&$ 0.482$&$0.058$&$0.486$&$ 8.340$&$0.000$&$0.000$&$0.059$&$0.059$\tabularnewline
$0.022$&$ 0.001$&$0.076$&$0.076$&$ 0.009$&$0.948$&$0.950$&$0.076$&$0.077$\tabularnewline
$0.042$&$ 0.001$&$0.078$&$0.078$&$ 0.008$&$0.944$&$0.948$&$0.076$&$0.077$\tabularnewline
$0.052$&$-0.001$&$0.078$&$0.078$&$-0.010$&$0.940$&$0.948$&$0.076$&$0.078$\tabularnewline
$0.112$&$ 0.002$&$0.081$&$0.081$&$ 0.026$&$0.930$&$0.946$&$0.076$&$0.080$\tabularnewline
$0.122$&$ 0.001$&$0.080$&$0.080$&$ 0.018$&$0.936$&$0.953$&$0.076$&$0.081$\tabularnewline
$0.252$&$ 0.002$&$0.088$&$0.088$&$ 0.026$&$0.912$&$0.949$&$0.076$&$0.088$\tabularnewline
$0.262$&$ 0.001$&$0.087$&$0.087$&$ 0.008$&$0.908$&$0.952$&$0.076$&$0.088$\tabularnewline
$0.504$&$-0.001$&$0.109$&$0.109$&$-0.013$&$0.827$&$0.950$&$0.076$&$0.108$\tabularnewline
$0.514$&$ 0.001$&$0.108$&$0.108$&$ 0.012$&$0.832$&$0.953$&$0.076$&$0.109$\tabularnewline
$0.524$&$ 0.000$&$0.110$&$0.110$&$ 0.003$&$0.825$&$0.948$&$0.076$&$0.110$\tabularnewline
$0.534$&$-0.004$&$0.110$&$0.110$&$-0.033$&$0.818$&$0.950$&$0.076$&$0.111$\tabularnewline
$0.544$&$ 0.001$&$0.111$&$0.111$&$ 0.012$&$0.819$&$0.949$&$0.076$&$0.112$\tabularnewline
$0.554$&$-0.001$&$0.111$&$0.111$&$-0.006$&$0.817$&$0.956$&$0.076$&$0.114$\tabularnewline
\hline
\end{tabular}

}\qquad\subfloat
[Model 6: Heteroskedastic $v_i$]{
\begin{tabular}{rrrrrrrrr}
\hline\hline
\multicolumn{1}{c}{$K/n$}&\multicolumn{1}{c}{$\mathsf{Bias}$}&\multicolumn{1}{c}{$\mathsf{SD}$}&\multicolumn{1}{c}{$\mathsf{RMSE}$}&\multicolumn{1}{c}{$\frac{\mathsf{Bias}}{\mathsf{SD}}$}&\multicolumn{1}{c}{CI$_0$}&\multicolumn{1}{c}{CI$_1$}&\multicolumn{1}{c}{$\hat\sigma$}&\multicolumn{1}{c}{$s$}\tabularnewline
\hline
$0.012$&$ 0.483$&$0.062$&$0.487$&$ 7.811$&$0.000$&$0.000$&$0.059$&$0.060$\tabularnewline
$0.022$&$ 0.001$&$0.077$&$0.077$&$ 0.011$&$0.945$&$0.948$&$0.076$&$0.077$\tabularnewline
$0.042$&$ 0.001$&$0.077$&$0.077$&$ 0.011$&$0.945$&$0.951$&$0.076$&$0.078$\tabularnewline
$0.052$&$-0.001$&$0.079$&$0.079$&$-0.009$&$0.941$&$0.948$&$0.077$&$0.079$\tabularnewline
$0.112$&$ 0.000$&$0.082$&$0.082$&$ 0.001$&$0.938$&$0.954$&$0.077$&$0.082$\tabularnewline
$0.122$&$ 0.004$&$0.080$&$0.080$&$ 0.046$&$0.942$&$0.955$&$0.077$&$0.082$\tabularnewline
$0.252$&$ 0.000$&$0.092$&$0.092$&$ 0.002$&$0.904$&$0.946$&$0.078$&$0.090$\tabularnewline
$0.262$&$ 0.002$&$0.089$&$0.089$&$ 0.026$&$0.910$&$0.957$&$0.078$&$0.091$\tabularnewline
$0.504$&$-0.001$&$0.117$&$0.117$&$-0.005$&$0.826$&$0.946$&$0.080$&$0.114$\tabularnewline
$0.514$&$-0.002$&$0.116$&$0.116$&$-0.017$&$0.828$&$0.951$&$0.081$&$0.116$\tabularnewline
$0.524$&$ 0.000$&$0.118$&$0.118$&$ 0.003$&$0.821$&$0.945$&$0.081$&$0.117$\tabularnewline
$0.534$&$ 0.001$&$0.118$&$0.118$&$ 0.010$&$0.815$&$0.953$&$0.081$&$0.119$\tabularnewline
$0.544$&$ 0.000$&$0.119$&$0.119$&$-0.003$&$0.816$&$0.952$&$0.081$&$0.120$\tabularnewline
$0.554$&$ 0.000$&$0.125$&$0.125$&$ 0.001$&$0.797$&$0.943$&$0.081$&$0.121$\tabularnewline
\hline
\end{tabular}

}\\
\footnotesize\raggedright Notes:\newline(i) columns $\mathsf{Bias}%
$, $\mathsf{SD}$, $\mathsf{RMSE}$ and $\frac{\mathsf{Bias}}{\mathsf{SD}%
}%
$ report, respectively, average bias, average standard deviation, root mean square error, and average standarized bias of the estimator $\hat
{\beta}$ across simulations;\newline
(ii) columns CI$_0$ and CI$_1$ report empirical coverage for homoskedastic-consistent confidence intervals, respectively, without and with degrees of freedom correction;\newline
(iii) columns $\hat\sigma
$ and $s$ report the average across simulations of the standard errors estimators, respectively, without and with degrees of freedom correction.
\end{table}%

\end{document}